\documentclass[11pt]{article}
\pdfoutput=1

\usepackage{mathrsfs}
\usepackage{amsmath,amssymb}
\usepackage{bm}
\usepackage{natbib}
\usepackage[usenames]{color}
\usepackage{amsthm}
\usepackage{algorithm}
\usepackage{algpseudocode}
\usepackage{multirow} 
\usepackage{enumitem}
\usepackage{caption}
\usepackage{subcaption}
\usepackage{enumitem}
\usepackage{dsfont}
\usepackage{mathtools}
\usepackage{caption}
\usepackage{booktabs}

\usepackage[colorlinks,
linkcolor=red,
anchorcolor=blue,
citecolor=blue
]{hyperref}

\usepackage{mylatexstyle}
\usepackage{float}  
\usepackage{setspace}
\usepackage[left=1in, right=1in, top=1in, bottom=1in]{geometry}

\setlist[itemize]{itemsep=0pt, topsep=2pt}
\setlist[enumerate]{itemsep=-2pt, topsep=2pt}

\usepackage{xcolor}

\ifdefined\final
\usepackage[disable]{todonotes}
\else
\usepackage[textsize=tiny]{todonotes}
\fi
\allowdisplaybreaks

\title{Fundamental Limits of Honest Inference
with Heterogeneous Transfer Sources}

\author
{
Xuran Meng\thanks{Department of Biostatistics, University of Michigan; e-mail: {\tt xuranm@umich.edu}}
	~~~~~~~and~~~~~~~
Yi Li\thanks{Department of Biostatistics, University of Michigan;
  e-mail: {\tt yili@umich.edu}}
}

\date{}

\DeclareMathOperator{\TV}{TV}
\DeclareMathOperator{\KL}{KL}

\newcommand{\Pn}{\mathsf P_n}

\newcommand{\E}{\mathbb E}

\begin{document}

\maketitle

\begin{abstract}
    Transfer learning can improve precision by borrowing information from external sources, but how much improvement is possible while maintaining valid inference? We address this by studying honest confidence intervals for a target parameter using multiple sources with heterogeneous prespecified discrepancy bounds. Under local information and admissibility conditions, we establish a universal lower bound on interval length for procedures maintaining uniform coverage. A scalar information equation characterizes the bound, revealing how precision depends jointly on source sample sizes and discrepancies. Large source samples need not yield greater precision gains when their discrepancies are large. When discrepancies are unrestricted, honest intervals cannot improve upon the target-only rate, even when all sources are unbiased. For the Gaussian location model, we calculate optimal source weights and construct honest intervals attaining the lower bound order, establishing sharpness. We extend the lower bound theory to semiparametric models through parametric submodels. For the Cox proportional hazards model, we verify the information conditions under which the lower bound applies, including under right censoring.
\end{abstract}

\section{Introduction}
\label{sec:introduction}

Auxiliary data can improve statistical precision when a target population is difficult or expensive to sample. However, the gains from borrowing information depend on the discrepancies between source and target populations. A large source sample may estimate its own parameters precisely, yet even a small discrepancy can undermine inference for the target. In particular, a discrepancy negligible relative to target sampling variability may exceed the width of an interval constructed primarily from source data. An honest procedure must maintain coverage across all configurations allowed by the discrepancy bounds. This raises a fundamental question: how much can external data improve statistical precision while maintaining uniform coverage, and what could limit these gains?

We address this question by studying honest confidence intervals for a target parameter using multiple sources with heterogeneous, prespecified discrepancy bounds, including infinite bounds representing unrestricted  discrepancies. Under local information and admissibility conditions, we establish a universal lower bound on interval length. A nonlinear scalar information equation characterizes this bound, with each source's contribution determined jointly by its sample size, discrepancy bound, and the configuration of the remaining sources. Consequently, a large source with a loose discrepancy bound may contribute less to honest inference than a smaller source with a tighter bound. When discrepancies are unrestricted, honest intervals cannot improve upon the target-only sampling rate, even when all sources are unbiased.

Our theory accommodates increasing parameter dimensions and
numbers of sources under uniform information and admissibility
conditions. In a Gaussian location model, we establish
sharpness of the lower bound by constructing honest intervals
attaining its order and deriving an optimal
source-weighting rule. This rule assigns zero weight to sources
whose discrepancy bounds exceed a common threshold.
We extend the theory to semiparametric models
through parametric submodels. For Cox proportional hazards
models \citep{cox1972regression}, we verify the information
condition under right censoring, complementing
the asymptotic analysis of \citet{andersen1982cox}
and high-dimensional estimation results of
\citet{huang2013oracle}. Counterexamples
illustrate the importance of the interior discrepancy
condition and show that the information condition
does not guarantee a matching upper bound.

Our work connects transfer learning with classical theory
on honest inference. In domain adaptation,
\citet{ben2010theory} investigate distributional differences,
while \citet{mansour2008domain} consider multiple sources.
For high-dimensional regression, \citet{li2022transfer}
establish estimation and prediction guarantees using auxiliary
samples, while \citet{tian2023transfer} study generalized
linear models, including source detection and coordinatewise
inference. Classical modulus-of-continuity and minimax theory
establish limits on honest inference
\citep{donoho1994statistical,low1997nonparametric,
cai2004adaptation}, with results for linear functionals
\citep{armstrong2018optimal} and high-dimensional regression
\citep{cai2017confidence}. Most closely related,
\citet{chen2025minimax} study minimax rates and adaptation
when combining experimental and observational data with
unknown confounding bias. Their confidence interval lower
bounds concern a scalar target, while their extension to
multiple biased estimators focuses on estimation.
In contrast, we establish lower bounds for honest inference
with multiple heterogeneous sources.

 Beyond these general minimax principles, our contribution is an explicit characterization of the limits of honest inference with heterogeneous, potentially biased sources. The resulting nonlinear information equation reveals how the effective contributions of different sources are coupled through a common threshold, providing a theoretical foundation for optimal source utilization. The Gaussian construction demonstrates the sharpness of this characterization through an explicit optimal weighting rule.

\section{Problem Formulation}
\label{sec:formulation}

\subsection{Target and Source Populations}

Let $p=p_n$ denote the parameter dimension and $K=K_n$ the number of sources, both indexed by $n$.
Consider a target cohort with sample size $n_0$ and distribution parameterized by $\beta\in\mathbb R^p$, and source cohorts indexed by $s=1,\ldots,K$, each with sample size $n_s$ and distribution parameterized by $\beta+\Delta_s$.
% Write distributions $R_\beta=R_{\beta,\nu_0}$ and $S_{s,\beta}=S_{s,\beta,\nu_s}$ for simplicity, as nuisance parameters $\nu_s$ are  fixed. When $\nu_s$ are absent, this model also reduces to the parametric case. 
We write the joint distribution of all the cohorts as 
\[
P_{n,\vartheta}
=R_\beta^{n_0}\times
S_{1,\beta+\Delta_1}^{n_1}\times\cdots\times
S_{K,\beta+\Delta_K}^{n_K},
\qquad
\vartheta=(\beta,\Delta_1,\ldots,\Delta_K).
\]
Here, \(R_\beta^{n_0}\) and \(S_{s,\beta+\Delta_s}^{n_s}\) denote the respective product distributions of \(n_0\) and \(n_s\) independent observations.  Without loss of generality, we consider inference for a fixed
coordinate $\beta_j$, where $1\le j\le p$, and let $e_j$
denote its coordinate vector. The same lower bounds apply
to the Euclidean diameter of honest confidence regions
for any fixed dimensional subvector.
% \ylcm{mention that the results hold (in what sense) for joint inference of a finite dimensional coefficients.}

Let \(\mathcal B_s\subseteq\mathbb R^p\) denote the parameter domain for cohort \(s\), with \(s=0\) representing the target. These domains satisfy the regularity conditions of the corresponding models. For known discrepancy bounds \(d_s\in[0,\infty]\), define the admissible parameter class
\[
\mathcal A_n=
\left\{
\vartheta:
\beta\in\mathcal B_0,\quad
\beta+\Delta_s\in\mathcal B_s,\quad
|\Delta_{s,j}|\le d_s,\quad 1\le s\le K
\right\}.
\]
 Setting \(d_s=\infty\) allows unrestricted or unknown discrepancies for source \(s\). We require confidence intervals to maintain uniform coverage over the admissible parameter class, as formalized below.
\begin{definition}
\label{def:honest_interval}
Fix \(0<\alpha<1/2\). An interval rule \(C_{n,j}\), constructed from the observed data, is honest at confidence level \(1-\alpha\) over \(\mathcal A_n\) if
\[
\inf_{\vartheta\in\mathcal A_n}
P_{n,\vartheta}\{\beta_j\in C_{n,j}\}
\ge 1-\alpha.
\]
\end{definition}
 The interval procedure may depend on the known parameter class and sample sizes, but not on the unknown true parameters. It may also use independent randomization, which is included in the calculation of probabilities and expectations. Let \(|C_{n,j}|\) denote the interval length, possibly infinite.

To evaluate interval length, let
\(\vartheta_n^0=(\beta^0,\Delta_1^0,\ldots,\Delta_K^0)\in\mathcal A_n\)
denote the true configuration, and write
\[
P_n=P_{n,\vartheta_n^0},\qquad
\beta_s^0=\beta^0+\Delta_s^0,\qquad
\delta_s^0=\Delta_{s,j}^0,\quad 1\le s\le K.
\]
We establish lower bounds on
\(P_n\{|C_{n,j}|\ge c\rho_n\}\) and
\(E_{P_n}|C_{n,j}|\) at true configurations satisfying
\[
|\delta_s^0|\le c_*d_s,\qquad 1\le s\le K,
\]
for a fixed \(c_*\in(0,1)\), subject to the stated information and admissibility conditions. The definition of $\rho_n$ is given in Section~\ref{sec:main}. This restriction places the true discrepancies away from the boundaries of positive finite discrepancy bounds. Section~\ref{sec:counterexamples} shows that an interval can be honest over \(\mathcal A_n\) yet have length \(o_p(\rho_n)\) at a boundary configuration.
All asymptotic statements concern a sequence of such experiments, typically with \(n_0\to\infty\), although the underlying information inequalities hold in finite samples.

\subsection{Comparison Configurations}

 To establish the lower bound, we compare the true configuration with an alternative in which the target coordinate is shifted by \(t\in\mathbb R\) and the corresponding coordinate in source \(s\) by \(u_s\). Let \(\mathbf u=(u_1,\ldots,u_K)\in\mathbb R^K\), and define
\[
\begin{aligned}
\Delta_s'&=\Delta_s^0+(u_s-t)e_j,\\
\vartheta_n(t,\mathbf u)
&=(\beta^0+te_j,\Delta_1',\ldots,\Delta_K').
\end{aligned}
\]
Under this alternative, the target parameter is \(\beta^0+te_j\), while the parameter in source \(s\) is \(\beta_s^0+u_se_j\). The corresponding discrepancy is \(\delta_s'=\delta_s^0+u_s-t\), wherein \(u_s=0\) leaves the source distribution unchanged. We require the alternative configuration to satisfy the following admissibility condition.

 \begin{assumption}[Admissibility of the comparison sequence]
\label{ass:comparison_admissibility}
Let \((t_n,\mathbf u_n)\), where \(\mathbf u_n=(u_{n,1},\ldots,u_{n,K})\), denote the comparison shifts used in the lower bound argument, satisfying
\(
|t_n|\vee\max_{1\le s\le K}|u_{n,s}|=O(1).
\)
We  assume that, for all sufficiently large \(n\), it holds that
\(
\beta^0+t_ne_j\in\mathcal B_0,  
\beta_s^0+u_{n,s}e_j\in\mathcal B_s,
   1\le s\le K.
\)
\end{assumption}

The allowable comparison shifts depend on the parameter domains. We verify Assumption~\ref{ass:comparison_admissibility} for the models considered in Section~\ref{sec:cases}. For such a comparison, write
\(
Q_n(t_n,\mathbf u_n)
=P_{n,\vartheta_n(t_n,\mathbf u_n)}.
\)
Membership of the alternative configuration in \(\mathcal A_n\) additionally requires
\[
|\delta_s^0+u_{n,s}-t_n|\le d_s,
\, 1\le s\le K.
\]
Under the interior condition \(|\delta_s^0|\le c_*d_s\), this holds whenever
\(
|u_{n,s}-t_n|\le(1-c_*)d_s,
\, 1\le s\le K.
\)
Thus, when both conditions are satisfied, honesty gives that, for all sufficiently large \(n\), 
\begin{equation}
\label{eq>=1-alpha}
Q_n(t_n,\mathbf u_n)
\{\beta_j^0+t_n\in C_{n,j}\}
\ge1-\alpha.
\end{equation}

The same interval procedure is applied under \(P_n\) and \(Q_n(t_n,\mathbf u_n)\), without knowledge of the true configuration or comparison shifts. Write \(x_+=\max(x,0)\). For probability laws \(P\) and \(Q\), define
\[
\operatorname{TV}(P,Q)=\sup_A|P(A)-Q(A)|,\qquad
D_{\mathrm{KL}}(P\Vert Q)
=E_P\log\left(\frac{dP}{dQ}\right).
\]
The lower bound follows by combining the coverage requirements under \(P_n\) and \(Q_n(t_n,\mathbf u_n)\) with a bound on the distance between these distributions. The key step is to select comparison shifts that separate the target coordinates by order \(\rho_n\) while keeping the two distributions sufficiently close. %The main thing in the proof is to construct the values of $t_n$ and $\ub_n$, .

\section{Main results}
\label{sec:main}
 We impose the following local information condition to control the
statistical distance between the true and alternative configurations.

\begin{assumption}[Uniform local information]
\label{assump:information}
There exists a constant $\kappa>0$, independent of $n,p,K$, and $s$,
such that, for all sufficiently large $n$,
\begin{equation*}
\begin{split}
D_{\KL}(R_{\beta^0}\Vert R_{\beta^0+a e_j})
&\le \kappa a^2,\\
D_{\KL}(S_{s,\beta_s^0}\Vert S_{s,\beta_s^0+a e_j})
&\le \kappa a^2,\qquad 1\le s\le K,
\end{split}
\end{equation*}
for all $|a|\le n_0^{-1/2}$ whenever the shifted parameter
belongs to the corresponding parameter domain.
\end{assumption}
% We establish a coordinatewise lower bound on the length of honest
% confidence intervals. At the actual configuration $\vartheta_n^0$,
% define $\rho_n$ as the positive solution of

% Since $F_n$ is continuous and strictly increasing on $[0,\infty)$,
% with $F_n(0)=0$ and $F_n(r)\to\infty$ as $r\to\infty$, this solution
% exists and is unique. Moreover, $ 0<\rho_n\le n_0^{-1/2}$.  
% The lower bound requires a local information condition and admissibility of the comparison configuration used in its proof.
 \begin{theorem}[Discrepancy-dependent lower bound]
\label{thm:multi}
Suppose Assumptions~\ref{ass:comparison_admissibility}
and \ref{assump:information} hold, and that the true
configuration satisfies
\(
|\delta_s^0|\le c_*d_s,   1\le s\le K,
\)
for some fixed $c_*\in(0,1)$.
Let $\rho_n$ be the unique positive solution to
\begin{equation}
F_n(\rho_n)
:=n_0\rho_n^2+
\sum_{s=1}^K n_s(\rho_n-d_s)_+^2=1.
\label{eq:root}
\end{equation}
Then there exists a constant $c>0$, depending only on
$\kappa$, $c_*$, and $\alpha$, such that every interval
$C_{n,j}$ that is honest at confidence level $1-\alpha$
over $\cA_n$ satisfies
\begin{equation*}
\begin{split}
\liminf_n \Pn\{|C_{n,j}|\ge c\rho_n\}
\ge \frac{1-2\alpha}{2},\quad
\liminf_n \frac{\E_{\Pn}|C_{n,j}|}{\rho_n}
\ge \frac{c(1-2\alpha)}{2}.
\end{split}
\end{equation*}
\end{theorem}

The scale $\rho_n$ characterizes how sample sizes and source
discrepancies jointly determine the lower bound, allowing both
the parameter dimension and the number of sources to grow.
Its implicit characterization admits a simpler order-equivalent
expression in terms of ordered discrepancy bounds and cumulative
source sample sizes.

 %Section~\ref{sec:cox} verifies the semi parametric model conditions for Cox proportional hazards models.  %The following proposition gives the order of \(\rho_n\).
\begin{proposition}
\label{prop:rho_n-sameorder}
Order the discrepancy bounds as
$d_{(1)}\le\cdots\le d_{(K)}$, with corresponding
sample sizes $n_{(1)},\ldots,n_{(K)}$, and set $d_{(0)}=0$.
Then the solution $\rho_n$ of \eqref{eq:root} satisfies
\begin{equation}
\rho_n\le\eta_n\le2\rho_n,
\qquad
\eta_n=
\min_{0\le m\le K}
\left\{
d_{(m)}+
\left(n_0+\sum_{r=1}^m n_{(r)}\right)^{-1/2}
\right\}.
\label{eq:cutoff-scale}
\end{equation}
\end{proposition}
 
When source discrepancies are unrestricted or unknown,
$d_s=\infty$ for every source, and
Proposition~\ref{prop:rho_n-sameorder} gives
$\rho_n=n_0^{-1/2}$.
By Theorem~\ref{thm:multi}, no honest transfer learning
procedure can achieve interval length $o_p(n_0^{-1/2})$,
even with arbitrarily large source samples.
Also, a matching upper bound is not guaranteed in general,
as illustrated by the counterexamples in
Section~\ref{sec:counterexamples}.

\section{Tightness, model examples and semiparametric extensions}
\label{sec:cases} 
\subsection{Tightness of the lower bound}
\label{sec:single}

We construct honest intervals attaining the lower bound order
in a Gaussian location model. Consider independent observations
$Y_{si}\sim N_p(\beta^0+\Delta_s^0,I_p)$,
$i=1,\ldots,n_s$, $s=0,\ldots,K$, with $\Delta_0^0=0$, 
unrestricted parameter domains and
$|\Delta_{s,j}^0|\le d_s$.
Since the parameter domains are unrestricted, Assumption~\ref{ass:comparison_admissibility} holds automatically.
$D_{\KL}\{N_p(\beta^0,I_p)\Vert
N_p(\beta^0+ae_j,I_p)\}=a^2/2$ verifies
Assumption~\ref{assump:information} with $\kappa=1/2$.
Set $d_0=0$ and $z=z_{1-\alpha/2}$.
For deterministic weights $w_s\ge0$ with $\sum_{s=0}^K w_s=1$,
define
\[
 \widehat\beta_j(w)=\sum_{s=0}^K w_s\bar Y_{s,j},
 \qquad
 r(w)=z(\sum_{s=0}^K{w_s^2}/{n_s})^{1/2}
       +\sum_{s=1}^K w_sd_s.
\]
Sources with $d_s=\infty$ receive zero weight and are omitted
from these sums.
The estimator is Gaussian with variance
$\sum_s w_s^2/n_s$ and absolute bias at most
$\sum_s w_sd_s$, so
$[\widehat\beta_j(w)-r(w),\widehat\beta_j(w)+r(w)]$
is honest over $\cA_n$. The minimum of $r(w)$ is attained by $w_s^*=
 \frac{n_s(r_*-d_s)_+}
 {\sum_{\ell=0}^K n_\ell(r_*-d_\ell)_+}$, and  $r_*=r(w^*)$ is the solution of
\begin{align}
 n_0r_*^2+\sum_{s=1}^K n_s(r_*-d_s)_+^2=z^2.
 \label{eq:optimal-radius}
\end{align}
Thus a source receives positive weight precisely when
$d_s<r_*$. The resulting interval
$C^{\mathrm G}_{n,j}
=[\widehat\beta_j(w^*)-r_*,\widehat\beta_j(w^*)+r_*]$
has deterministic length $2r_*$. From Equation~\eqref{eq:optimal-radius}, we have  $2\min\{1,z\}\rho_n
 \le |C^{\mathrm G}_{n,j}|
 \le 2\max\{1,z\}\rho_n$.  The lower bound order in Theorem~\ref{thm:multi} is thus
attained in the Gaussian location model.

% When $d_{1}\lesssim n_1^{-1/2}$, the lower bound scale is
% $\rho_n\asymp n_1^{-1/2}$. If
% $d_{1}\lesssim n_1^{-1/2}$, the expanded source interval
% attains this order.
% When $n_1^{-1/2}\ll d_{1}\ll n_0^{-1/2}$, the scale is
% $\rho_n\asymp d_{1}$. If $d_{1}\asymp d_{1}$,
% the added discrepancy allowance dominates the source sampling
% radius, giving length of order $d_{1}$.
% Finally, when $d_{1}\gtrsim n_0^{-1/2}$,
% $\rho_n\asymp n_0^{-1/2}$, and the target-only interval attains
% this order. Formula~\eqref{eq:single-scale} also covers the
% boundaries between these regimes. 
% Thus, under the stated conditions, the construction achieves
% $|C_{n,j}|=O_{p}(\rho_n)$.
% Under the assumptions of Theorem~\ref{thm:multi}, every honest
% interval satisfies
% \[
%  \Pn\{|C_{n,j}|\ge c\rho_n\}\ge\frac{1-2\alpha}{2},
% \]
% so $|C_{n,j}|=o_{p}(\rho_n)$ is impossible.
% The lower bound is therefore sharp in order in these regimes.
\subsection{Some Examples}
\label{sec:counterexamples}

We construct two examples. The first example demonstrates the necessity of the
interior discrepancy condition. Consider scalar Gaussian cohorts with unit variance,
two sources with $n_1=n_2=N$, and fixed bounds $d_1=d_2=d>1$.
Suppose $n_0\to\infty$ and $N/n_0\to\infty$, so that
$\rho_n=n_0^{-1/2}$ eventually. Put
$r=z_{1-\alpha/2}/\sqrt N$,
$L=\bar Y_1-d-r$, $U=\bar Y_2+d+r$, and
$C=[L,\max\{L,U\}]$.
Since $\beta_1^0-d\le\beta^0\le\beta_2^0+d$, Gaussian tail bounds
give coverage at least $1-\alpha$ throughout the class
$|\delta_1^0|,|\delta_2^0|\le d$.
At the boundary configuration $(\delta_1^0,\delta_2^0)=(d,-d)$,
however,
\[
 |C|=[\bar Y_2-\bar Y_1+2d+2r]_+
 =O_p(N^{-1/2})=o_p(\rho_n).
\]
Thus the pointwise lower bound cannot extend to every boundary configuration.

The second example shows that the information condition alone also does not guarantee a matching
upper bound. Suppose every cohort follows
\[
 R_{\beta}
 =N_2\bigl((\beta_1+\beta_2, 
           \varepsilon_n\beta_2)^\T,I_2\bigr),
 \qquad \varepsilon_n>0,\quad \varepsilon_n\to0,
\]
with unrestricted parameter domains, $j=1$, and $d_s=0$. The true configurations are $\beta^0=(\beta_1^0,\beta_2^0)$ and $\vartheta^0=(\beta^0,0,...,0)$. 
Under this case a matching coordinatewise KL lower bound holds.
Write $M=\sum_{s=0}^K n_s$, giving $\rho_n=M^{-1/2}$.
The configurations with common cohort parameters $(0,0)^\T$
and $(t,-t)^\T$ both belong to the coverage class.
Their joint laws $P$ and $Q$ satisfy
$D_{\KL}(P\Vert Q)=M\varepsilon_n^2t^2/2$.
Taking $t=(1-2\alpha)/(\varepsilon_n\sqrt M)$,
every interval $C_{n,j}$ honest at confidence level $1-\alpha$ over $\cA_n$ satisfies
\[
 \E_{n,\vartheta^0}|C_{n,j}|
 \ge
 \frac{(1-2\alpha)^2}{2\varepsilon_n}\rho_n.
\]
The ratio to $\rho_n$ therefore diverges. A matching upper bound may thus require additional conditions.
\subsection{Semiparametric Extension: The Cox Model}
\label{sec:balls}
\label{sec:cox} 
The lower bound extends to semiparametric models through parametric submodels satisfying the stated assumptions. Within each cohort, we fix nuisance parameters at their true values across comparisons. These parameters may differ across cohorts and need not be known to the procedure. Any interval honest over the full semiparametric model must also satisfy the coverage requirement on these submodels.
 The Cox proportional hazards model  provides a semiparametric example. 
For cohort $s=0,\ldots,K$, suppose
\(
X_s\sim N_p(0,\Sigma_s)\), \(
\lambda_s(v\mid X_s;\beta)
=\lambda_{0s}(v)e^{X_s^\T\beta}\), \(
 \Lambda_{0s}(\infty)=\infty,
\)
where $\Lambda_{0s}(v)=\int_0^v\lambda_{0s}(r) dr$,
and the baseline hazard is nonnegative and locally integrable.
% The baseline hazard and covariate distribution are fixed within
% each comparison. 
Let $R^{\mathrm{full}}_{s,\beta}$ denote the joint law of $(T_s,X_s)$. 
Consider the parameter domains
$\cB_s=\{\beta:\|\beta\|_\infty\le M\}$,
for some $M>0$. The comparison admissibility condition holds if
\(
\max_{0\le s\le K}\|\beta_s^0\|_\infty\le M-\varepsilon
\)
for some  $\varepsilon>0$, provided $n_0\to\infty$
and the  shifts satisfy
\(
|t_n|\vee\max_{1\le s\le K}|u_{n,s}|
\le n_0^{-1/2}.
\)
Hence, each comparison changes  smaller than
$\varepsilon$ for all large $n$. In comparisons, nuisance parameters  are fixed,  and  Theorem~\ref{thm:multi} thus applies on the resulting parametric submodels.

\begin{proposition}[Cox information bound]
\label{prop:cox}
For any parameter shift $h\in\RR^p$, $D_{\KL}(R^{\mathrm{full}}_{s,\beta}\Vert
 R^{\mathrm{full}}_{s,\beta+h})
 =\exp(h^\T\Sigma_sh/2)-1$.
If $\sup_{n, 0\le s\le K_n}\Sigma_{s,jj}\le\Lambda<\infty$
for some $\Lambda>0$, then
\begin{equation}
 D_{\KL}(R^{\mathrm{full}}_{s,\beta}\Vert
 R^{\mathrm{full}}_{s,\beta+a e_j})
 \le\kappa a^2,\qquad |a|\le1,
 \qquad \kappa=\frac{\Lambda}{2}e^{\Lambda/2}.
 \label{eq:cox-quadratic}
\end{equation}
\end{proposition}
By Proposition~\ref{prop:cox}, a uniform bound on $\Sigma_{s,jj}$ suffices for
Assumption~\ref{assump:information}, with a constant independent
of the coefficients, baseline hazards, dimension and number of
sources. The proof is in Appendix~\ref{app:cox}.

 For right censored observations
$O_s=(T_s\wedge C_s,\mathbf1\{T_s\le C_s\},X_s)$,
suppose $C_s\perp T_s\mid X_s$ and the conditional censoring
law is unchanged across comparison configurations.
By the data processing inequality,
\begin{equation}
\begin{split}
D_{\KL}(R^{\mathrm{obs}}_{s,\beta}\Vert
R^{\mathrm{obs}}_{s,\beta+a e_j})
\le
D_{\KL}(R^{\mathrm{full}}_{s,\beta}\Vert
R^{\mathrm{full}}_{s,\beta+a e_j}) \le \kappa a^2,\qquad |a|\le1.
\end{split}
\label{eq:censoring}
\end{equation}
where $R^{\mathrm{obs}}_{s,\beta}$ denotes the law of $O_s$.
Thus, Assumption~\ref{assump:information} also holds
under right censoring.
% The baseline hazards and covariance matrices may be unknown: honesty over the larger model implies honesty on the submodel where they are fixed.

\section{Experiments}
\label{sec:experiments}
 We compare the Gaussian intervals in Section~\ref{sec:single}
with target-only intervals, both constructed to achieve at least
$95\%$ coverage uniformly over $\cA_n$. The corresponding
inverted tests therefore have type I error at most $0.05$,
although actual coverage may vary across configurations.

We simulate independent Gaussian sample means
$\bar Y_{s,j}\sim N(\delta_s,n_s^{-1})$,
$s=0,\ldots,K$, with $\beta_j=\delta_0=0$.
We consider $K\in\{5,20\}$ and
$n_0=n\in\{100,400,1600,6400\}$.
The source sample sizes are $n_s=n^{a_s}$,
where $a_s$ increases linearly from $1.5$ to $3$
across the $K$ sources. For $s<K$, the discrepancy
bounds are $d_s=n^{-b_s}$, where $b_s$ decreases
linearly from $1$ to $0.2$ across the first $K-1$
sources, and $d_K=10$.

The true discrepancies are $\delta_s=(-1)^s c d_s$,
with $c\in\{0.5,0.8\}$, satisfying the interior
discrepancy condition. The optimal transfer weights
depend only on the known discrepancy bounds and
sample sizes. The target-only interval is
$[\bar Y_{0,j}-z_{0.975}/\sqrt n,
  \bar Y_{0,j}+z_{0.975}/\sqrt n]$.
Table~\ref{tab:supplied} reports the results.
% Each setting uses $20000$ independent replications at  confidence level $0.95$. We report empirical coverage, average interval length, the length divided by $\rho_n$, and the ratio of transfer to target-only length.  The actual discrepancies are used only to generate the data and assumed to be unknown. 

\begin{table}[t]
\caption{Interval lengths and empirical type I error from $20,000$
Gaussian replications. Both procedures control type I error
at 0.05. Here, $L$ and $L_0$ denote transfer and target-only
interval lengths, $R_n=L/\rho_n$, and
$m_\rho=\#\{s:d_s<\rho_n\}$.
The last block uses the first $k$ sources with $K=5$ and $n=1600$.}
\label{tab:supplied}
\centering
\small
\setlength{\tabcolsep}{4pt}
\renewcommand{\arraystretch}{1.1}
\begin{tabular}{@{}rrrrrrcccc@{}}
\toprule
& & & & & & \multicolumn{4}{c}{Empirical type I error}\\
\cmidrule(l){7-10}
& & & & & &
\multicolumn{2}{c}{Transfer} &
\multicolumn{2}{c}{Target-only}\\
\cmidrule(lr){7-8}\cmidrule(l){9-10}
$n$ & $m_\rho$ & $\rho_n$ & $L$ & $R_n$ & $L/L_0$
& $c=0.5$ & $c=0.8$ & $c=0.5$ & $c=0.8$\\
\midrule
\multicolumn{10}{@{}l}{$K=5$}\\
\addlinespace[2pt]
100
& 2 & $3.79\times10^{-2}$ & $1.04\times10^{-1}$ & 2.75 & .266
& .001 & .003 & .048 & .050\\
400
& 2 & $1.31\times10^{-2}$ & $3.48\times10^{-2}$ & 2.67 & .178
& .000 & .003 & .050 & .051\\
1600
& 2 & $4.51\times10^{-3}$ & $1.18\times10^{-2}$ & 2.61 & .120
& .000 & .003 & .047 & .051\\
6400
& 1 & $1.54\times10^{-3}$ & $4.01\times10^{-3}$ & 2.60 & .082
& .000 & .001 & .050 & .050\\
% 25600
% & 1 & $5.31\times10^{-4}$ & $1.38\times10^{-3}$ & 2.60 & .056
% & .000 & .001 & .048 & .049\\
\midrule
\multicolumn{10}{@{}l}{$K=20$}\\
\addlinespace[2pt]
100
& 5 & $2.57\times10^{-2}$ & $6.90\times10^{-2}$ & 2.68 & .176
& .000 & .000 & .049 & .048\\
400
& 5 & $7.88\times10^{-3}$ & $2.15\times10^{-2}$ & 2.73 & .110
& .000 & .000 & .049 & .050\\
1600
& 5 & $2.42\times10^{-3}$ & $6.68\times10^{-3}$ & 2.76 & .068
& .000 & .000 & .051 & .052\\
6400
& 4 & $7.40\times10^{-4}$ & $2.06\times10^{-3}$ & 2.78 & .042
& .000 & .001 & .050 & .053\\
% 25600
% & 4 & $2.27\times10^{-4}$ & $6.36\times10^{-4}$ & 2.81 & .026
% & .000 & .001 & .052 & .047\\
\midrule
\multicolumn{10}{@{}l}{$K=5$, $n=1600$, first $k$ sources included ($k=0$: target-only)}\\
\addlinespace[3pt]
$k$ & $m_\rho$ & $\rho_n$ & $L$ & $R_n$ & $L/L_0$
& $c=0.5$ & $c=0.8$ & $c=0.5$ & $c=0.8$\\
\cmidrule{1-10}
0
& 0 & $2.50\times10^{-2}$ & $9.80\times10^{-2}$ & 3.92 & 1.000
& .047 & .051 & .047 & .051\\
1
& 1 & $4.51\times10^{-3}$ & $1.65\times10^{-2}$ & 3.66 & .169
& .035 & .036 & .047 & .051\\
2
& 2 & $4.51\times10^{-3}$ & $1.18\times10^{-2}$ & 2.61 & .120
& .000 & .003 & .047 & .051\\
3
& 2 & $4.51\times10^{-3}$ & $1.18\times10^{-2}$ & 2.61 & .120
& .000 & .003 & .047 & .051\\
4
& 2 & $4.51\times10^{-3}$ & $1.18\times10^{-2}$ & 2.61 & .120
& .000 & .003 & .047 & .051\\
5
& 2 & $4.51\times10^{-3}$ & $1.18\times10^{-2}$ & 2.61 & .120
& .000 & .003 & .047 & .051\\
\bottomrule
\end{tabular}
\end{table}

As shown in Table~\ref{tab:supplied}, $R_n$ ranges from $2.60$ to $2.78$, within the theoretical interval $[2,3.92]$. These nearly constant ratios show the length of order $\rho_n$ achieved by the Gaussian construction.
By Theorem~\ref{thm:multi}, no honest procedure can achieve length $o_p(\rho_n)$.
The construction is therefore optimal in order, and the experiments illustrate its attainment of that rate.
The transfer to target-only length ratios are also well below one. For example, with $K=20$ and $n=6400$, the transfer interval is only about $4.2\%$ compared to the target-only interval. Empirical transfer type I error ranges from $0.000$ to $0.003$, while target-only type I error remains around $0.05$. Therefore, even heterogeneous source data with known discrepancy bounds  could potentially provide  gains in precision while controlling type I error. The last panel of Table~\ref{tab:supplied} considers first $k$ sources included with $K=5$ and $n=1600$.
With all five sources included, $\rho_n\approx4.51\times10^{-3}$,
and only the first two sources satisfy $d_s<\rho_n$. We can see that adding the first source reduces the interval length from
$9.80\times10^{-2}$ to $1.65\times10^{-2}$, and adding the
second reduces it further to $1.18\times10^{-2}$.
The third to fifth sources have much larger sample sizes,
but their bounds also exceed the optimal radius $r_*$,
so they receive zero weight. 
This illustrates the threshold structure of the lower bound
equation~\eqref{eq:root} and the optimal Gaussian weighting
rule~\eqref{eq:optimal-radius}.
%This matches our theoretical equation~\eqref{eq:root} that only sources with discrepancy bounds below the  threshold $\rho_n$ increase the information.

%\section{Discussion}
%\label{sec:discussion}
%We establish a lower bound on honest confidence interval length that depends on the sample sizes and discrepancies of all sources. Under the stated assumptions, no procedure with uniform coverage over admissible coverage class can achieve length $o_{p}(\rho_n)$. We also construct an honest interval that matches the lower bound $\Theta_p(\rho_n)$  in the Gaussian model in Section~\ref{sec:single}. Thus the lower bound shows a tight limit on the general honest procedure. Experiments further validate our theory.

\appendix

\section{Proofs}
\label{app:proofs}

\subsection{Proof of Theorem~\ref{thm:multi}}
\label{app:multi}
  We construct the comparison shifts
\(
c=(1-c_*)\min\left\{\frac12,
\frac{1-2\alpha}{\sqrt{2\kappa}}\right\},
\,
t_n=c\rho_n,\,
u_{n,s}=(t_n-(1-c_*)d_s)_+,
\, 1\le s\le K.
\)
Write $\ub_n=(u_{n,1},\ldots,u_{n,K})^\T$.
Since $\rho_n\le n_0^{-1/2}$, we have
\(
0\le u_{n,s}\le t_n\le c n_0^{-1/2}.
\)
Moreover, the interior discrepancy condition implies
\[
|\delta'_{n,s}|
\le |\delta_s|+|u_{n,s}-t_n|
\le c_*d_s+(1-c_*)d_s=d_s.
\]
Thus, Assumption~\ref{ass:comparison_admissibility}
ensures that the comparison configuration
$\vartheta_n(t_n,\ub_n)$ belongs to $\cA_n$
for all sufficiently large $n$, while the shifts
satisfy the range required by
Assumption~\ref{assump:information}.
Let $Q_n=P_{n,\vartheta_n(t_n,\ub_n)}$.
Finally, since $c\le 1-c_*$, we have 
\(
u_{n,s}
=(c\rho_n-(1-c_*)d_s)_+
\le c(\rho_n-d_s)_+.
\)
 
By independence and Assumption~\ref{assump:information}, we have
\vspace{-0.1cm}
\begin{align*}
D_{\KL}(\Pn\Vert Q_n)
&=n_0D_{\KL}(R_{\beta_0}\Vert R_{\beta_0+t_ne_j})
+\sum_{s=1}^K n_s
D_{\KL}(S_{s,\beta_s}\Vert S_{s,\beta_s+u_{n,s}e_j})\\
&\le \kappa\left(n_0t_n^2+
\sum_{s=1}^K n_su_{n,s}^2\right)
\le \kappa c^2\left\{
n_0\rho_n^2+
\sum_{s=1}^K n_s(\rho_n-d_s)_+^2
\right\}
=\kappa c^2,
\end{align*}
where the last equality follows from \eqref{eq:root}.
Pinsker's inequality therefore gives
\[
\TV(\Pn,Q_n)
\le \sqrt{\frac{D_{\KL}(\Pn\Vert Q_n)}{2}}
\le c\sqrt{\frac{\kappa}{2}}
\le \frac{1-2\alpha}{2}.
\]

Let $A_0=\{\beta_{0j}\in C_{n,j}\}$ and
$A_1=\{\beta_{0j}+t_n\in C_{n,j}\}$.
Since both configurations belong to $\cA_n$,
honesty implies
\(
\Pn(A_0)\ge 1-\alpha,\,
Q_n(A_1)\ge 1-\alpha.
\)
By the total variation bound,
\[
\Pn(A_1)\ge Q_n(A_1)-\TV(\Pn,Q_n)
\ge 1-\alpha-\TV(\Pn,Q_n).
\]
On $A_0\cap A_1$, the interval contains two points
separated by $t_n=c\rho_n$. Hence
\begin{align*}
\Pn\{|C_{n,j}|\ge c\rho_n\}
&\ge \Pn(A_0\cap A_1)\\
&\ge \Pn(A_0)+\Pn(A_1)-1\\
&\ge 1-2\alpha-\TV(\Pn,Q_n)
\ge \frac{1-2\alpha}{2}.
\end{align*}
Consequently,
\[
\frac{\E_{\Pn}|C_{n,j}|}{\rho_n}
\ge c\,\Pn\{|C_{n,j}|\ge c\rho_n\}
\ge \frac{c(1-2\alpha)}{2}.
\]
This completes the proof of Theorem~\ref{thm:multi}.
\hfill$\square$

\subsection{Proof of Proposition~\ref{prop:rho_n-sameorder}}
\label{app:cutoff}

For $0\le m\le K$, define $N_m=n_0+\sum_{r=1}^m n_{(r)}
 $, $
 q_m=d_{(m)}+N_m^{-1/2},$ 
with an empty sum when $m=0$.
Thus $\eta_n=\min_{0\le m\le K}q_m$.
For each $r\le m$, $q_m-d_{(r)}\ge N_m^{-1/2}$, $
  q_m\ge N_m^{-1/2}$.  
It follows that $F_n(q_m)
 \ge
 n_0q_m^2+
 \sum_{r=1}^m n_{(r)}(q_m-d_{(r)})^2
 \ge N_mN_m^{-1}=1$.  
This also holds for $m=0$, when $q_0=n_0^{-1/2}$.
Since $F_n$ is strictly increasing and $F_n(\rho_n)=1$,
we have $\rho_n\le q_m$ for every $m$, and hence $\rho_n\le\eta_n$.

For the reverse inequality, let $m_*$ be the number of sources
with $d_{s}<\rho_n$.
If $m_*=0$, it is clear that $\rho_n=n_0^{-1/2}=\eta_n$. 
If $m_*>0$, the root equation gives $1
 =n_0\rho_n^2+
 \sum_{r=1}^{m_*}n_{(r)}
       (\rho_n-d_{(r)})^2
 \le N_{m_*}\rho_n^2$.  
Thus $N_{m_*}^{-1/2}\le\rho_n$.
Since $d_{(m_*)}<\rho_n$, $ \eta_n
 \le q_{m_*}
 =d_{(m_*)}+N_{m_*}^{-1/2}
 \le2\rho_n$. 
These two inequalities prove \eqref{eq:cutoff-scale}.
\hfill$\square$

\subsection{Proof of Proposition~\ref{prop:cox}}
\label{app:cox}

Suppress the cohort index $s$. The conditional density is
$f_\beta(v\mid x)=\lambda_0(v)e^{x^\T\beta}
\exp\{-\Lambda_0(v)e^{x^\T\beta}\}$ for $v\ge0$.
Since $\Lambda_0(\infty)=\infty$,
$V=\Lambda_0(T)e^{X^\T\beta}$ has a unit exponential distribution
conditional on $X$, so $\E_\beta(V\mid X)=1$.
With the covariate law and baseline hazard fixed, the log
likelihood ratio is 
\begin{align*}
 \log\frac{f_\beta(T\mid X)}{f_{\beta+h}(T\mid X)}
 =-X^\T h+
   \Lambda_0(T)e^{X^\T\beta}(e^{X^\T h}-1)=-X^\T h+V(e^{X^\T h}-1).
\end{align*}
% Taking expectation conditional on $X$ replaces $V$ by its
% conditional mean, giving $e^{X^\T h}-1-X^\T h$.
Since $X\sim N_p(0,\Sigma)$, we have
$\E(X^\T h)=0$ and
$\E e^{X^\T h}=\exp(h^\T\Sigma h/2)$. Hence $D_{\KL}(R^{\mathrm{full}}_\beta\Vert
         R^{\mathrm{full}}_{\beta+h})
 =\E\{e^{X^\T h}-1-X^\T h\}
 =\exp(h^\T\Sigma h/2)-1$.
For $h=ae_j$ with $|a|\le1$, put $x=\Sigma_{jj}a^2/2$.
The variance bound gives $0\le x\le\Lambda/2$, and
$e^x-1\le xe^x\le(\Lambda/2)e^{\Lambda/2}a^2$.
Thus \eqref{eq:cox-quadratic} holds with
$\kappa=(\Lambda/2)e^{\Lambda/2}$.
The same constant applies to every cohort and parameter value
with uniformly bounded  $\Sigma_{s,jj}$.

For censored observations, let $H(dc\mid x)$ denote the common
conditional law of $C$ given $X=x$.
Conditional independence implies that the joint law of $(T,X,C)$
is $R^{\mathrm{full}}_\beta(dt,dx)H(dc\mid x)$.
The factor $H(dc\mid x)$ cancels from the likelihood ratio.
The observed data are obtained through the common map
$(T,X,C)\mapsto(T\wedge C,\mathbf1\{T\le C\},X)$.
The data processing inequality then gives $D_{\KL}(R^{\mathrm{obs}}_\beta\Vert
         R^{\mathrm{obs}}_{\beta+h})
 \le
 D_{\KL}(R^{\mathrm{full}}_\beta\Vert
         R^{\mathrm{full}}_{\beta+h})$. 
Taking $h=ae_j$ and applying the quadratic bound proves
\eqref{eq:censoring}.
\hfill$\square$

\bibliography{ref}

\begin{thebibliography}{13}
\expandafter\ifx\csname natexlab\endcsname\relax\def\natexlab#1{#1}\fi
\expandafter\ifx\csname url\endcsname\relax
  \def\url#1{\texttt{#1}}\fi
\expandafter\ifx\csname urlprefix\endcsname\relax\def\urlprefix{URL }\fi

\bibitem[{Andersen and Gill(1982)}]{andersen1982cox}
\textsc{Andersen, P.~K.} and \textsc{Gill, R.~D.} (1982).
\newblock Cox's regression model for counting processes: a large sample study.
\newblock \textit{The Annals of Statistics}  1100--1120.

\bibitem[{Armstrong and Koles{\'a}r(2018)}]{armstrong2018optimal}
\textsc{Armstrong, T.~B.} and \textsc{Koles{\'a}r, M.} (2018).
\newblock Optimal inference in a class of regression models.
\newblock \textit{Econometrica} \textbf{86} 655--683.

\bibitem[{Ben-David et~al.(2010)Ben-David, Blitzer, Crammer, Kulesza, Pereira and Vaughan}]{ben2010theory}
\textsc{Ben-David, S.}, \textsc{Blitzer, J.}, \textsc{Crammer, K.}, \textsc{Kulesza, A.}, \textsc{Pereira, F.} and \textsc{Vaughan, J.~W.} (2010).
\newblock A theory of learning from different domains.
\newblock \textit{Machine learning} \textbf{79} 151--175.

\bibitem[{Cai and Guo(2017)}]{cai2017confidence}
\textsc{Cai, T.~T.} and \textsc{Guo, Z.} (2017).
\newblock Confidence intervals for high-dimensional linear regression: Minimax rates and adaptivity.
\newblock \textit{The Annals of Statistics} \textbf{45} 615--646.

\bibitem[{Cai and Low(2004)}]{cai2004adaptation}
\textsc{Cai, T.~T.} and \textsc{Low, M.~G.} (2004).
\newblock An adaptation theory for nonparametric confidence intervals.
\newblock \textit{The Annals of Statistics} \textbf{32} 1805--1840.

\bibitem[{Chen et~al.(2025)Chen, Li, Zhang and Ye}]{chen2025minimax}
\textsc{Chen, S.}, \textsc{Li, S.}, \textsc{Zhang, B.} and \textsc{Ye, T.} (2025).
\newblock Minimax rates and adaptivity in combining experimental and observational data.
\newblock \textit{Journal of Causal Inference} \textbf{13} 0024.

\bibitem[{Cox(1972)}]{cox1972regression}
\textsc{Cox, D.~R.} (1972).
\newblock Regression models and life-tables.
\newblock \textit{Journal of the royal statistical society: Series B (methodological)} \textbf{34} 187--202.

\bibitem[{Donoho(1994)}]{donoho1994statistical}
\textsc{Donoho, D.~L.} (1994).
\newblock Statistical estimation and optimal recovery.
\newblock \textit{The Annals of Statistics} \textbf{22} 238--270.

\bibitem[{Huang et~al.(2013)Huang, Sun, Ying, Yu and Zhang}]{huang2013oracle}
\textsc{Huang, J.}, \textsc{Sun, T.}, \textsc{Ying, Z.}, \textsc{Yu, Y.} and \textsc{Zhang, C.-H.} (2013).
\newblock Oracle inequalities for the lasso in the cox model.
\newblock \textit{Annals of statistics} \textbf{41} 1142.

\bibitem[{Li et~al.(2022)Li, Cai and Li}]{li2022transfer}
\textsc{Li, S.}, \textsc{Cai, T.~T.} and \textsc{Li, H.} (2022).
\newblock Transfer learning for high-dimensional linear regression: Prediction, estimation and minimax optimality.
\newblock \textit{Journal of the Royal Statistical Society Series B: Statistical Methodology} \textbf{84} 149--173.

\bibitem[{Low(1997)}]{low1997nonparametric}
\textsc{Low, M.~G.} (1997).
\newblock On nonparametric confidence intervals.
\newblock \textit{The Annals of Statistics} \textbf{25} 2547--2554.

\bibitem[{Mansour et~al.(2008)Mansour, Mohri and Rostamizadeh}]{mansour2008domain}
\textsc{Mansour, Y.}, \textsc{Mohri, M.} and \textsc{Rostamizadeh, A.} (2008).
\newblock Domain adaptation with multiple sources.
\newblock \textit{Advances in neural information processing systems} \textbf{21}.

\bibitem[{Tian and Feng(2023)}]{tian2023transfer}
\textsc{Tian, Y.} and \textsc{Feng, Y.} (2023).
\newblock Transfer learning under high-dimensional generalized linear models.
\newblock \textit{Journal of the American Statistical Association} \textbf{118} 2684--2697.

\end{thebibliography}
\bibliographystyle{ims}

\end{document}